\documentclass[12pt,reqno]{amsproc}

\usepackage[utf8]{inputenc}
\usepackage[all,cmtip]{xy}    

\usepackage{enumerate}
\usepackage{amsfonts, mathpazo}
\usepackage{bbm}
\usepackage{amssymb}
\usepackage{amsmath}
\usepackage{mathrsfs}
\usepackage{amsthm}

\usepackage{changes}           
\usepackage[margin=.9in]{geometry}
\usepackage{graphicx}
\usepackage{subfig}
\usepackage{hyperref}
\hypersetup{
    colorlinks=true,
    linkcolor=blue,
    filecolor=blue,      
    urlcolor=red,
    citecolor=red,
    }
\usepackage[capitalise]{cleveref}
\usepackage{todonotes}
\usepackage{orcidlink}

\usepackage{tikz}
\usepackage{tikz-cd}  
\usepackage{caption}
\usepackage{subcaption}

\newcommand{\define}[1]{{\bf \boldmath{#1}}}
\newcommand{\set}[1]{\left\{{#1}\right\}}
\newcommand{\card}[1]{\left\lvert{#1}\right\rvert}

\def \R {{\mathbb R}}
\def \Z {{\mathbb Z}}

\def \eps {{\varepsilon}}

\def \cl {{\mathrm{\bf cl}}}

\def \argmin {{\mathrm{argmin}}}

\def \sd {{\mathrm{Sd}}}
\def \st {\mathrm{\bf st}}
\def \lk {{\mathrm{\bf lk}}}

\def \cone {{\mathrm{cone}}}
\def \Con {{\textrm{\bf Con}}}
\def \HG {{\textrm{H}}}
\def \ex {{\textrm{\bf ex}}}
\def \cM {\mathscr{M}}

\newcommand{\halo}{{\mathbf{h}}}
\newcommand{\ahalo}{\widehat{\halo}}
\newcommand{\res}{\mathbf{sh}}
\newcommand{\U}[1]{{\mathbf{U}{#1}}}
\newcommand{\coll}{\searrow}

\newcommand{\scoll}{\raisebox{1ex}{\rotatebox{-45}{\ensuremath{\rightsquigarrow}}}~}

\theoremstyle{plain}
\newtheorem{theorem}{Theorem}[section]
\newtheorem*{theorem*}{Theorem}
\newtheorem{lemma}[theorem]{Lemma}
\newtheorem{proposition}[theorem]{Proposition}

\theoremstyle{definition}
\newtheorem{definition}[theorem]{Definition}

\newtheorem*{example*}{Example}
\newtheorem{remark}[theorem]{Remark}

\title{Stratified Morse Theory for Cell Complexes}
\author{Vidit Nanda and Francesca Tombari}

\begin{document}

\begin{abstract}
We develop a version of discrete Morse theory for finite regular CW complexes equipped with an auxiliary stratification. The key construction is the \emph{halo} of a cell, which contains all those faces in the boundary that enter closed sublevelsets precisely when the threshold reaches that cell's value. The complement of this halo in the boundary, called the \emph{shadow}, is always a subcomplex. A stratified discrete Morse function requires Forman's conditions on each stratum together with the requirement that closures of paired cells admit filtered collapses onto their shadows. We establish fundamental Morse lemmas: filtered collapses across regular intervals, and controlled attachments at critical values. For functions satisfying only the stratum-wise Forman condition, we construct an upper envelope on the barycentric subdivision whose local Morse data decomposes into horizontal and vertical components. This yields a simplicial analogue of the standard tangential-normal splitting of local Morse data in the sense of Goresky and MacPherson.

\end{abstract}

\maketitle

\section*{Introduction}

Morse theory provides one of the finest and most successful mechanisms for translating local analytic information into global geometric insight. Beginning with Morse's original work on counting geodesics \cite{morse}, the core insight -- that the topology of sublevelsets changes only at critical points, and does so in a controlled manner -- has also found spectacular applications when extended to equivariant \cite{atiyahbott}, piecewise-linear \cite{bb}, symplectic \cite{floer}, dynamical \cite{conley}, and physical \cite{witten} contexts. In each case, the fundamental theorems survive: regular intervals produce no topological change, while critical values yield handle attachments whose nature can be determined by local data of the underlying function near critical points. We are concerned here with two such extensions --- the first of these is the combinatorial adaptation developed by Forman for finite CW complexes~\cite{forman1998}, and the second is Goresky-MacPherson's Morse theory for stratified spaces \cite{GoreMPer1988}. 

In Forman's theory, a discrete Morse function $f$ assigns real values to cells of a CW complex $X$ such that each cell $\sigma$ has at most one exceptional coface $\tau > \sigma$ with $f(\tau) \leq f(\sigma)$, or at most one exceptional face $\eta < \sigma$ with $f(\eta) \geq f(\sigma)$. Cells with no exceptional neighbours are deemed critical, while the remaining cells assemble into free-face pairs that can be collapsed away. This discrete avatar of Morse theory has also found substantial applications across diverse fields, including topological combinatorics \cite{babson}, commutative algebra \cite{welker}, algebraic topology \cite{salvetti, braid}, and geometric group theory \cite{adiprasito, dmcog}. It has also been modified in several directions --- there are filtered \cite{dmpers}, equivariant \cite{Freij2009, dmcog}, sheafy \cite{dmsheaf}, noncompact \cite{kukiela} and even 2-categorical refinements \cite{flow1, flow2}. Our goal here is to develop a \emph{stratified} discrete Morse theory that combines insights from Goresky and MacPherson's pioneering work~\cite{GoreMPer1988} while retaining the combinatorial essence and flexibility of Forman's.

Therefore, the setting throughout is that of a finite regular CW complex $X$ equipped with a filtration
\[
X = X_n \supset X_{n-1} \supset \cdots \supset X_0 \supset X_{-1} = \varnothing
\]
such that the successive differences $X_i - X_{i-1}$, whose connected components are called {\em strata}, satisfy the frontier axiom (see Section \ref{ssec:strat}). This induces a partial order on the strata. Besides this order, we emphasise at the outset that  we impose no manifold structure, no regularity conditions \`{a} la Whitney or Thom--Mather, and no requirement that strata admit conical neighbourhoods. We forsake these hypotheses deliberately, for two reasons: first, we do not need them. And second, imposing them would prevent our theory from recovering Forman's results for the trivial stratification $X \supset \varnothing$.

\subsection*{The Challenge} 

Our primary target in this work is the pair of {\bf fundamental Morse lemmas}. The first of these asserts that nothing of topological importance changes across regular intervals --- sublevelsets are homotopy equivalent, and indeed, collapse onto one another in a filtration-preserving manner. The second one quantifies the change across a critical value: the sublevelset grows by a controlled handle attachment. Before we can establish (or even formulate) discrete analogues of these results, we must confront two obstacles. 

The first of these is conceptual: {\em how should we even define a discrete stratified Morse function} $f:X \to \R$? One might na\"ively require that the restriction $f|_S:S \to \R$ to every stratum $S \subset X$ satisfies the axioms of Forman. Such a strategy is antithetical to both the spirit and the letter of Goresky-MacPherson's setting, where a crucial additional \emph{normal nondegeneracy} requirement is imposed. The gradient of a stratified Morse function $f$ at a critical point must not annihilate any limiting tangent space of a higher stratum at that point. This nondegeneracy condition ensures that gradient flow does not become trapped at stratum boundaries. {\em What is the combinatorial equivalent of this condition}? We have no tangent planes, no gradients in the differential-geometric sense, and no obvious way to define or detect transversality.

\begin{figure}[ht]
\includegraphics[scale=.6]{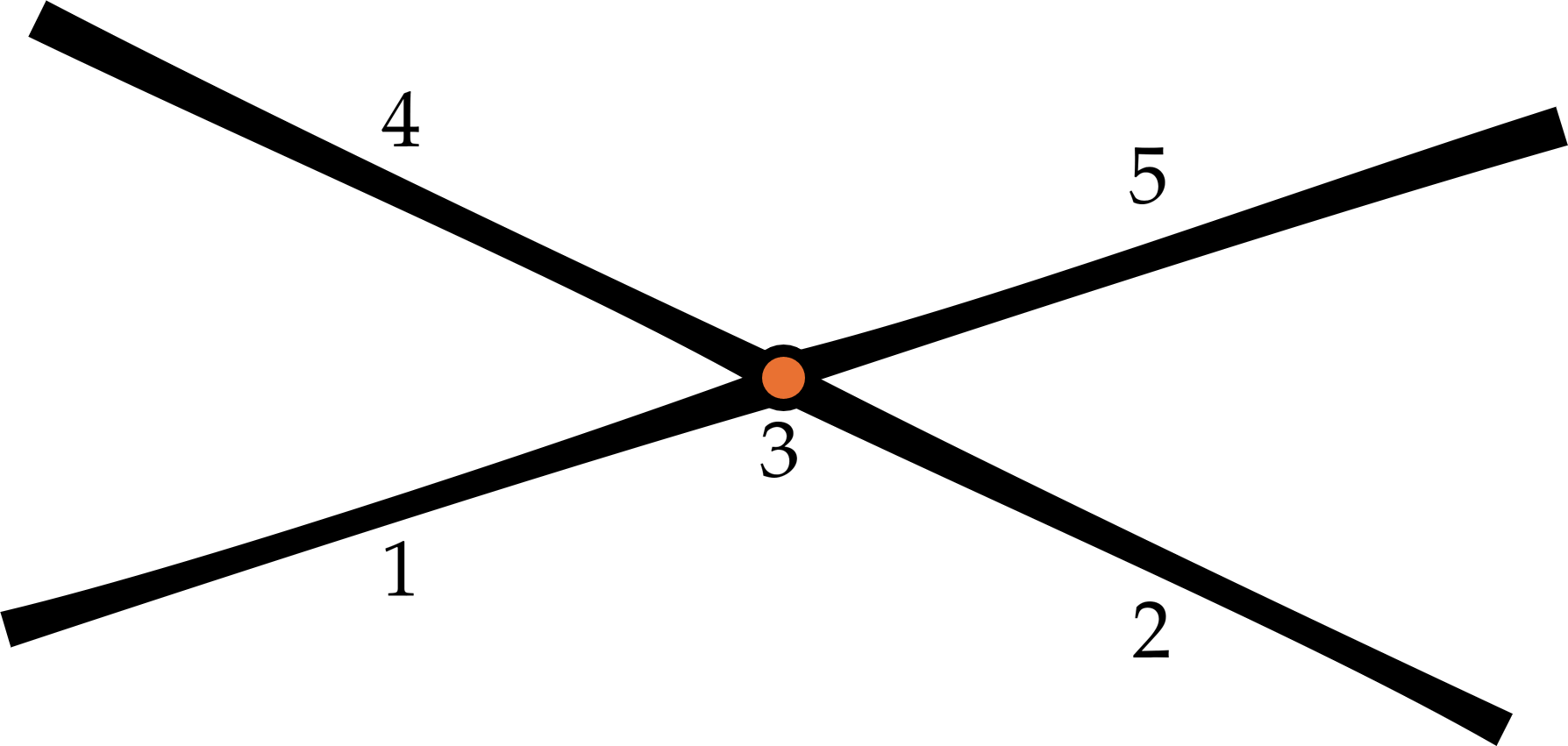}
\caption{A piece of a cell complex containing one vertex and four edges, with $f$-values indicated on the cells. Note that the closure of $f_{\leq 1}$ already contains the central vertex, even though that vertex has $f$-value $3$.}
\label{fig:closure}
\end{figure}

The second obstacle is technical, and we call it the \emph{closure problem}. In Goresky and MacPherson's setting, the sublevelset $f_{\leq c}$ is unambiguous --- it is simply the preimage of the half-infinite interval $(-\infty, c]$ under a continuous function. In sharp contrast, sublevelsets of the form $f_{\leq c}$ need not form subcomplexes of the ambient CW complex $X$ precisely because we allow cells $\sigma < \tau$ to satisfy $f(\sigma) \geq f(\tau)$. Following Forman, we are therefore compelled to work with closures: the relevant object is $\cl(f_{\leq c})$, the smallest subcomplex containing all cells with $f$-value at most $c$. The problem is that passing to the closure can introduce cells with arbitrarily large $f$-values --- a cell $\sigma$ with $f(\sigma) = 10000$ will appear in $\cl(f_{\leq 1})$ if $\sigma$ has a coface $\tau > \sigma$ with $f(\tau) \leq 1$ (see Figure \ref{fig:closure}). There is no parallel to this phenomenon in the continuous setting, one might as well compare a volcano to a butterfly.

\subsection*{This Paper} 

Here we resolve both obstacles by introducing the \emph{halo} $\halo(\sigma;f)$ of a given cell $\sigma \in X$ along an injective map $f:X \hookrightarrow \R$. This is defined as the collection of all boundary cells $\eta < \sigma$ such that $f(\eta) \geq f(\sigma)$ and $\sigma$ is the $f$-minimal coface of $\eta$. Informally, the halo consists of those faces that appear in the closed sublevelset $\cl(f_{\leq c})$ precisely when the threshold $c$ reaches $f(\sigma)$. The complement of the halo in the boundary of $\sigma$, which we call the \emph{shadow} of $\sigma$ and denote $\res(\sigma;f)$, forms a subcomplex of the closure $\cl(\sigma)$. The halo allows us to formulate our combinatorial nondegeneracy condition as follows. In addition to requiring our putative Morse function to be injective and satisfy the Forman axioms on each stratum, we further ask that for each cell pair $(\sigma > \tau)$ lying in the same stratum with $f(\sigma) < f(\tau)$, the smaller cell $\tau$ lies in $\halo(\sigma;f)$ and the closure $\cl(\sigma)$ admits a filtered collapse onto the shadow $\res(\sigma; f)$.

This new condition ensures that when a cell $\sigma$ is cancelled by its partner $\tau$ within a stratum, all the collateral faces in the halo can also be collapsed away without disrupting the stratification. The filtered collapse respects the filtration $X_\bullet$ in the sense that it restricts to a homotopy equivalence at each level. As a result, critical cells now come in two flavours: there is the standard notion, where a cell $\sigma$ lying in a stratum $S$ is \emph{critical} if it is unpaired within $S$. It is {\em s-critical} if it is critical and does not lie in the halo of any other cell. This distinction is genuinely new: a cell can be critical and contribute meaningfully to the topology of its stratum, and yet be invisible to the filtered topology of $X$ because it admits cofaces with smaller $f$-values lying in higher strata. Only s-critical cells affect the filtered homotopy type of sublevelsets.

\subsection*{Main Results}

Fix a stratified discrete Morse function $f:X \hookrightarrow \R$ and let $\sigma$ be a cell of $X$. Our first main result establishes that there is no material change in the sublevelsets of $f$ across intervals which contain no s-critical values. Here is an abridged version, see Theorem \ref{thm:a} for details.

\begin{theorem*}[A]
Let $[c, d] \subset \R$ be an interval whose preimage $f^{-1}\left([c,d]\right)$ is $\set{\sigma}$. If $\sigma$ is not s-critical, then $\cl(f_{\leq d})$ admits a filtration-preserving collapse onto $\cl(f_{\leq c})$.
\end{theorem*}

Let us emphasise once again that the collapse promised by this result might occur even if $\sigma$ is critical, for the reasons detailed above --- if $\sigma$ lies in the halo of another cell, then it will already be present in $\cl(f_{\leq c})$ and make no material contribution as we increase the threhold from $c$ to $d$. Having confirmed that sublevelsets can only change when we cross an s-critical value, it remains to produce {\bf local Morse data} by describing the precise handle attachment across such a value. This task becomes the purview of our next main result (see Theorem \ref{thm:b} below).

\begin{theorem*}[B]
Let $\sigma \in X$ be an s-critical cell of $f$ with $c := f(\sigma)$. For all sufficiently small $\eps > 0$, the sublevelset $\cl(f_{\leq c+\eps})$ is the union of $\cl(\sigma)$ with $\cl(f_{\leq c-\eps})$ along their common intersection $\res(\sigma; f)$.
\end{theorem*}

In other words, the {local Morse data} at an s-critical cell $\sigma$ is the pair $\res(\sigma; f) \hookrightarrow \cl(\sigma)$. We should specify that this data is manifestly different from Goresky and MacPherson's tangential $\times$ normal decomposition from \cite{GoreMPer1988}; their tangential data $T$ occurs along the stratum $S$ containing the critical point $p$ and is identical to the usual handle attachment data for the restriction $f|_S$. Their normal data $N$ is generated by the {\em descending link} of $S$, which resides in the strata lying strictly above $S$. This disparity is both a feature -- our attachment data is combinatorially explicit and requires no link computations -- and a limitation, in that we obtain less disentangled geometric information than the smooth theory provides. Our third and final main result is therefore a peace offering to the reader who would much prefer a discrete theory that faithfully recreates the continuous one's $T \times N$ product decomposition of local Morse data. 

In fact, here we consider a weaker setup. We still require $f:X \to \R$ to be injective and restrict to a discrete Morse function on each stratum, but no further halo-centric constraints are imposed.  We then pass to the barycentric subdivision of $X$; this is a simplicial complex $\sd(X)$ whose $k$-simplices are strictly ascending chains $\xi = [\sigma_0 < \cdots < \sigma_k]$ of cells in $X$. The subdivision canonically inherits its stratification from $X$, where each stratum $S \subset X$ induces a unique stratum $\sd(S) \subset \sd(X)$. The map $f$ extends to $\sd(X)$ via its {\em upper envelope} $\U{f}:\sd(X) \to \R$, given by
\[
\U{f}(\xi) = \max \set{f(\sigma_i) \mid 0 \leq i \leq k}.
\]  
The upshot is that sublevelsets of $\U{f}$ are honest subcomplexes and there is no need to take closures. We establish that -- for the sublevelsets of $\U{f}$ -- the local Morse data at a vertex $[\sigma]$ in $\sd(X)$ corresponding to each cell $\sigma \in X$ admits a join decomposition of the form $H_\sigma \star V_\sigma$, where $H = H_\sigma$ and $V = V_\sigma$ are subcomplexes called the {\em horizontal} and {\em vertical} parts of the data. The terminology is justified by the following observations: 
if $\sigma$ is critical and $S$ is the stratum containing it, then $H$ lies in the union of strata $\leq \sd(S)$ while $V$ lies in the union of strata $> \sd(S)$. 

\begin{theorem*}[C]
 Let $\sigma$ be a critical cell with $c := f(\sigma)$. For all sufficiently small $\eps > 0$, the local Morse data for $\U{f}_{\leq c-\eps} \hookrightarrow \U{f}_{\leq c+\eps}$ is a pair
\[
(H \star V) \hookrightarrow (\cone(H) \star V),
\]
where $H$ lies in the union of strata $\leq \sd(S)$ and $V$ lies in the union of strata $> \sd(S)$.
\end{theorem*}
\noindent (This is Theorem \ref{thm:tndecomp} below).

\begin{figure}[ht]
\includegraphics[scale=.6]{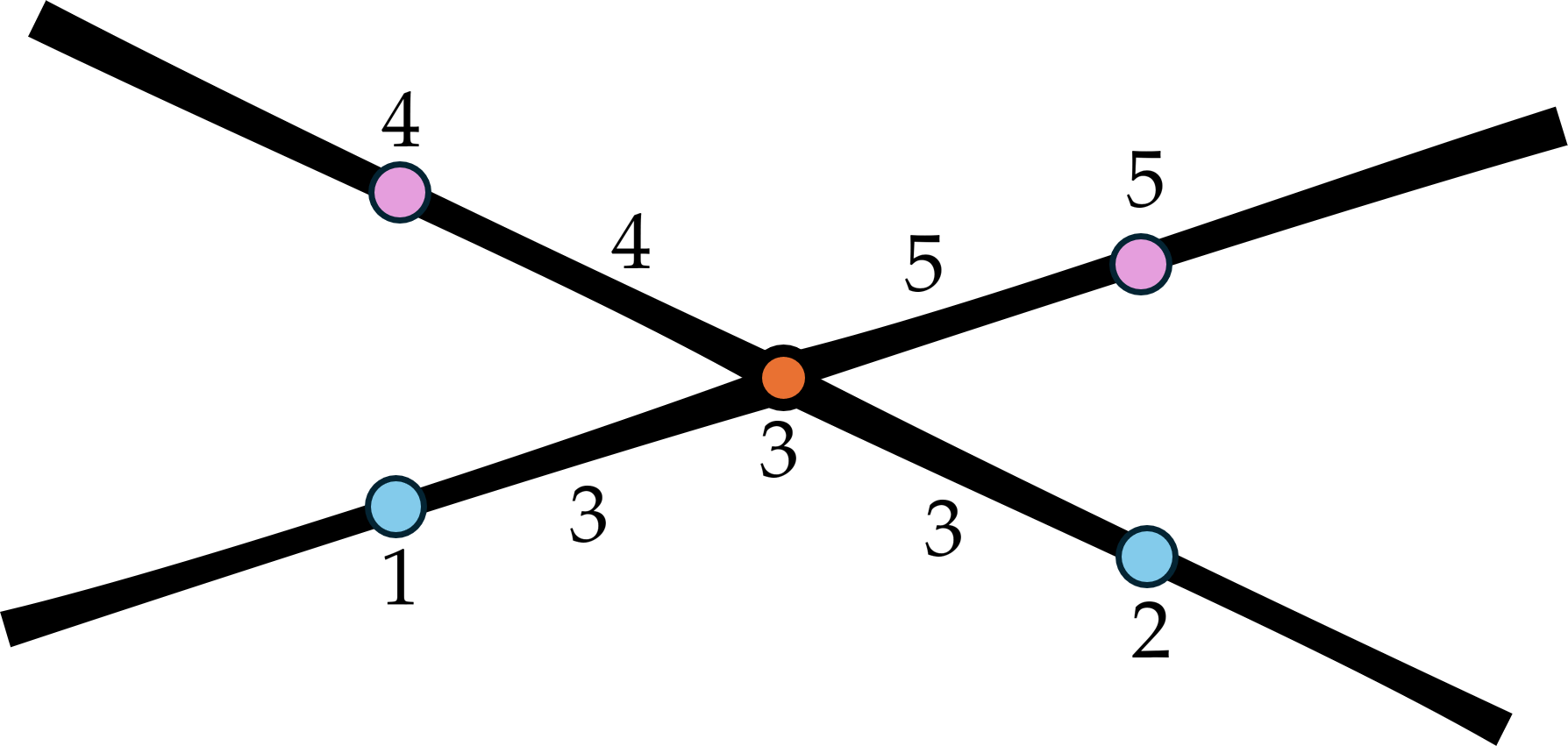}
\caption{A subdivided version of the Example from Figure \ref{fig:closure}, this time with $\U{f}$-values on simplices; here the central vertex does contribute to the change in topology because its lower link consists of two vertices. These lie at the barycenters of the original edges, and hence have $\U{f}$-values $1$ and $2$.}
\label{fig:subdivision}
\end{figure}

This decomposition is our simplicial analogue of the tangential-normal splitting. The horizontal part $H$ captures Morse data along the stratum, while the vertical part $V$ records the directions in which higher strata impinge on $\sigma$. 
In Figure \ref{fig:subdivision} we consider the barycentric subdivision and upper envelope of Figure \ref{fig:closure}.
The central vertex, which had no contribution in the original model of taking sublevelsets, now admits a lower link of the form $H \star V$, where $H$ is empty but $V$ has two connected components (these are barycenters of edges with $f$-values $1$ and $2$). Since this link is not contractible, crossing $\U{f} = 3$ is guaranteed to change sublevelset topology.

\subsection*{Related Work} 

This paper may be viewed as continuing three largely independent research programmes, all of which seek to generalise Forman's original theory from \cite{forman1998}. The first of these is the filtered variant due to Mischaikow and Nanda \cite{dmpers}. The main difference between that setting and ours is that in our work, we allow for the existence of cells $\sigma < \tau$ lying in different strata with $f(\sigma) > f(\tau)$ --- this is expressly prohibited in \cite[Def 4.1]{dmpers}. 

The second overlapping research direction is the more recent output of Knudson and Wang \cite{knudwang2022}, which also aims to build a stratified discrete Morse theory. Their functions are only required to restrict to Forman-style discrete Morse functions on each stratum, with no further requirements. As a consequence, they obtain weak Morse inequalities relating critical cell counts to Betti numbers, but are unable to extract local Morse data. Our Theorem (C) shows that even in their setting, one recovers a tangential-normal decomposition in the barycentric subdivision as long as one has injectivity.

Thirdly, we note the substantial literature on multivector fields, initiated by Mrozek and collaborators~\cite{mro2017, MrozekWanner2015Dynamics, LipKubMroWan201}, which provides a framework for discrete Conley index theory. A multivector is a nonempty, convex, connected subset of the face poset; a multivector field is a partition of a regular CW complex $X$ into multivectors. Stratifications satisfying the frontier axiom yield acyclic multivector fields, and the Conley--Morse spectral sequence assembles the homology of $X$ from the Conley indices of its strata. We discuss this connection briefly in an appendix, but emphasise that multivector field theory imposes no analogue of our halo collapsibility condition, and therefore does not yield Morse lemmas of the kind we establish here.

\subsection*{Outline} 

In \cref{sec:strat} we describe stratifications of regular CW complexes. Then \Cref{sec:halo} introduces halos and shadows, proving their basic properties. In \Cref{sec:stratmorse} we define stratified discrete Morse functions and the two notions of criticality. \Cref{sec:morselems} establishes Theorems (A) and (B). Finally, \Cref{sec:envelope} constructs the upper envelope on the barycentric subdivision and proves Theorem C. Appendix \ref{sec:conley} tersely summarises the relationship between our stratifications and acyclic multivector fields.

\section{Stratified Cell Complexes} \label{sec:strat}

Let $X$ be a regular CW complex. We denote by $(X,\leq)$ the \define{face poset} of $X$, which contains all the cells of $X$ ordered by the face relation\footnote{Explicitly, the relation $\sigma < \tau$ holds if and only if the boundary of $\tau$ contains $\sigma$.}. Given any nonempty subcollection $S$ of cells lying in $X$, the {\em closure} of $S$ is the subcomplex defined as a down-set
\[
\cl(S) := \set{\tau \in X \mid \tau \leq \sigma \text{ for some }\sigma \in S};
\]
and dually, the {\em (open) star} of $S$ is the up-set
\[
\st(S) := \set{\tau \in X \mid \tau \geq \sigma \text{ for some }\sigma \in S}.
\]
When $S = \set{\sigma}$ consists of a single cell, we  write $\cl(\sigma)$ and  $\st(\sigma)$ rather than $\cl(\set{\sigma})$ and $\st(\set{\sigma})$. We say that $S$ is {\em connected} if any two cells $\sigma,\tau$ in $S$ form endpoints of a finite zigzag:
\[
\sigma \geq \eta_0 \leq \eta_1 \geq \cdots \geq \eta_{k-1} \leq \eta_k \geq \tau,
\] with all $\eta_i \in S$. Finally, we denote by $\dim S$ the maximal dimension encountered among the constituent cells of $S$.

\subsection{Stratifications}\label{ssec:strat}

A {\em filtration} of a finite regular CW complex $X$ is a finite sequence of subcomplexes
\[
X = X_n \supset X_{n-1} \supset \cdots \supset X_1 \supset X_0 \supset  X_{-1}=\varnothing.
\]
A \define{stratification} of $X$ is a filtration
such that the successive differences $\Delta_i := (X_i - X_{i-1})$  satisfy the following {\em frontier axiom} --- given connected $S \subset \Delta_i$ and $T \subset \Delta_j$, if $\cl(S) \cap T$ is nonempty then $T \subset \cl(S)$. The connected components of $\Delta_i$ are called $i$-{\em strata} of the stratification $X_\bullet$. A standard example is the {\em skeletal} stratification, where the $i$-strata are precisely the $i$-dimensional cells; this is the finest possible stratification of $X$. At the other end of the spectrum, one has the coarsest (or  canonical, or minimal) stratification of $X$ into homology manifolds --- see \cite{locstrat}.

\begin{remark}
	As we have mentioned already in the Introduction, it is customary when defining stratified spaces to impose additional regularity conditions on the strata. For instance: one often requires $i$-strata to be (smooth, piecewise-linear, or homology) manifolds of dimension $i$ \cite{whitney, mather, weinberger, locstrat}. {\bf We make no such assumptions here}. Similarly, it is customary to assume that strata admit conical neighbourhoods in $X$. Namely, for each stratum $S$ there exists a stratified space $LS$, called the {\bf link} of $S$, such that any sufficiently small open neighbourhood $US \subset X$ of $S$ admits a fiber bundle structure $US \twoheadrightarrow S$ with fiber $\text{Cone}(LS)$. Our stratifications are less rigid in the sense that the strata are not required to admit well-defined links or homogeneous neighbourhoods.
\end{remark}

Nevertheless, we impose the frontier axiom in order to avail of two direct benefits. The first of these is the existence of a well-beloved and oft-used partial order on the strata.

\begin{proposition} \label{prop:frontier}
In any stratification $X_\bullet$ of $X$, the relation
\[
T \leq S \text{ if and only if } T \cap \cl(S) \neq \varnothing
\]
is a partial order on the set of strata. 
\end{proposition}
\begin{proof}
    Reflexivity of the proposed relation is immediate. 
    If $T\le S$ and $S\le R$, by definition of the relation and the frontier axiom, we have $T\subset \cl(S)$ and $S\subset \cl(R)$. 
    This implies $T\subset \cl(S)\subset \cl(R)$, and hence $T\le R$, showing transitivity.
    It remains to establish antisimmetry.
    Proceeding by contradiction, assume there exist strata $S \neq T$ such that both 
    $S\cap \cl(T)$ and $T\cap \cl(S)$ are nonempty. Let us assume that $S \subset \Delta_i$ and $T \subset \Delta_j$ where $\Delta_i = (X_i - X_{i-1})$ and similarly for $\Delta_j$. We now have containments
    \[
    T \subset \cl(S) \subset X_i, 
    \] where the first inclusion follows from the frontier axiom while the second follows because $X_i$ is a subcomplex, and hence closed in $X$. Thus, $T$ lies in $\Delta_j \cap X_i$, which forces $i \geq j$. The same argument gives $j \geq i$ if we exchange the roles of $S$ and $T$, whence $i=j$. Thus, both $S$ and $T$ must be connected components of the same $\Delta_i$. Finally, since $T \subset \cl(S)$, there exists cells $\tau \in T$ and $\sigma \in S$ with $\tau \leq \sigma$. This forces $S$ and $T$ to lie in the same connected component of $\Delta_i$, and hence $S=T$ as desired.
    \end{proof}

For the skeletal stratification, the frontier partial order evidently reduces to the familiar face relation between cells.
The second benefit of the frontier axiom is that it guarantees {\em convexity} of all strata, as described below.

\begin{proposition}\label{prop:convexstrata}
	Given any triple of cells $\sigma \leq \tau \leq \sigma'$ in a regular CW complex $X$, if both $\sigma$ and $\sigma'$ lie in the same stratum $S$ of a stratification $X_\bullet$, then so does $\tau$. 
\end{proposition}
\begin{proof}
	The only nontrivial case occurs when both inequalities among cells are strict, otherwise $\tau \in S$ is automatic. By construction, the strata of $X_\bullet$ partition the cells of $X$; we may therefore assume for the sake of contradiction that $\tau$ lies in a different stratum $T \neq S$. Now $T$ intersects $\cl(S)$ because $\tau < \sigma'$ whereas $S$ instersects $\cl(T)$ because $\sigma < \tau$. By the frontier axiom, we therefore get both $S \leq T$ and $T \leq S$, which immediately forces the desired contradiction $S = T$.
\end{proof}

\subsection{Filtered Collapses} \label{ssec:collapse} 

 We say that distinct cells $(\sigma,\tau)$ constitute a {\em free-face pair} in a regular CW complex $X$ whenever
\[
\st(\tau) = \set{\sigma,\tau}.
\]
In this case, the regularity of $X$ forces $\dim \tau = \dim \sigma-1$ and it is well-known (see eg \cite{cohen}) that removing both $\sigma$ and $\tau$ from $X$ produces a subcomplex which remains homotopy equivalent to $X$. This free-face pair removal operation is called an \define{elementary collapse}; we say that $X$ {\em collapses} to a subcomplex $Y$ if there is a finite sequence of intermediate subcomplexes $Z^i$
\[
X = Z^1 \supset Z^2 \supset \cdots  \supset Z^k = Y
\]
where each $Z^i$ has been obtained from the preceding $Z^{i-1}$ by performing a single elementary collapse. We denote this relationship between $X$ and $Y$ by $X \coll Y$.

\begin{remark}
	The reader unfamiliar with this machinery is warned that collapsing $X$ onto $Y$ is not simply a matter of finding a collection of disjoint free-face pairs in $X$. When a free-face pair $(\sigma,\tau)$ is removed from $X$ to produce the first subcomplex $Z^1$, the excision of $\sigma$ might create a new free-face pair $(\sigma',\tau')$ in $Z^1$ for some  $\tau' < \tau$. Crucially, this new $(\sigma',\tau')$ is {\em not} a free-face pair in $X$ because $\st(\tau') \supset \set{\sigma',\sigma,\tau'}$ has more than two elements.
\end{remark}

Let us assume now that $X$ is equipped with a stratification $X_\bullet$; every subcomplex $Y \subset X$ inherits this stratification via $Y_i := Y \cap X_i$ for all $i$. We now seek a more refined collapsing mechanism from $X$ to $Y$ which produces a {\em filtered} homotopy equivalence. Thus, the goal is to not only obtain a homotopy equivalence $\phi_i:X_i \stackrel{\sim}{\longrightarrow} Y_i$ at each level $i$, but also to require compatibility with the natural inclusion maps which relate adjacent levels. Explicitly, the diagram 
\[
\begin{tikzcd}
    X_i\ar[r, hook]\ar[d, "\phi_1"', "\sim"] & X_{i+1}\ar[d, "\sim"', "\phi_{i+1}"]\\
    Y_i\ar[r, hook] & Y_{i+1}
\end{tikzcd}
\]
should commute for all $i$. This task is accomplished by allowing the stratification data to constrain the class of acceptable elementary collapses --- we only permit the removal of a free-face pair $(\sigma,\tau)$ when both $\sigma$ and $\tau$ lie in the same stratum $S$. We will denote the existence of such a {\bf filtered collapse} as $X \scoll Y$ in order to distinguish it from the ordinary unfiltered collapses described above. 

\begin{remark} Assume that the regular CW subcomplex $Y \subset X$ has been obtained via a filtered elementary collapse. Denote by $(\sigma,\tau)$ the excised pair, let $S \subset (X_j - X_{j-1})$ be the stratum containing both cells. For each dimension $i$, we seek to relate the induced stratification $Y_i := Y \cap X_i$ of $Y$ to the original stratification of $X$ in order to make clear why a filtered homotopy equivalence $\phi_i:X_i \to Y_i$ exists. For this purpose, note first that $Y_i = X_i$ holds whenever $i < j$ simply because neither $\sigma$ nor $\tau$ are present in any such $X_j$, so we let $\phi_i$ be the identity for such $i$. Next, we note that $(\sigma,\tau)$ is a free-face pair in $X_i$ for all $i \geq j$ --- we have $\st(\tau) = \set{\sigma,\tau}$ in $X$, and since both $\sigma$ and $\tau$ lie in $X_j$ this pair remains free in $X_i$ for all $i \geq j$. It follows that one can choose the global homotopy equivalence $\phi:X \to Y$ to restrict to homotopy equivalences $X_i \to Y_i$ as expected, and hence the natural square commutes.
\end{remark}

\section{Halos and Shadows} \label{sec:halo}

Let $X$ be a finite regular CW-complex. Throughout this section, we consider a function $f:X \to \R$ that sends each cell $\sigma$ of $X$ to a real number $f(\sigma)$. The {\em lower star} of $\sigma$ along $f$ is
\begin{align}\label{eq:lowstar}
\st^-(\sigma;f) := \set{\tau \in X \mid \tau > \sigma \text{ and } f(\tau) \le f(\sigma)};
\end{align}
this is evidently a subset of $\st(\sigma)$, but -- unlike the usual star -- it never contains $\sigma$ itself. We will be particularly interested here in the case where $f$ is injective. For each collection of cells $S \subset X$, let $\argmin[f/S]$ denote the subset of all cells $\tau \in S$ which attain the minimal $f$-value in the image $f(S) \subset \R$, with the explicit understanding that $\argmin[f/\varnothing] = \varnothing$. For injective $f$ and nonempty $S$, the set $\argmin[f/S]$ is guaranteed to be a singleton.

\begin{definition}\label{def:halo}
The \define{halo} of a cell $\sigma \in X$ along $f$ is defined as 
\[
\halo(\sigma;f) := \set{\tau \in X \mid \sigma\in \argmin\left[f/\st^-(\tau;f)\right]},
\]
while the {\em augmented} halo $\ahalo(\sigma;f)$ is the (necessarily disjoint) union of $\halo(\sigma;f)$ with $\set{\sigma}$.
\end{definition}
 
Note that if $\tau$ lies in $\halo(\sigma;f)$ then we must have $\sigma \in \st^-(\tau;f)$. 
Thus, both $\halo(\sigma;f)$ and $\ahalo(\sigma;f)$ are always subsets (but generally not subcomplexes) of the closure $\cl(\sigma)$. For injective $f$, a cell $\tau$ lies in $\halo(\sigma;f)$ if and only if $\sigma$ is the unique cell in $\st^-(\tau;f)$ with minimal $f$-value. This forces the halos of distinct cells along injective $f$ to be disjoint.
\begin{proposition}\label{prop:halodisjoint}
	Let $f:X \hookrightarrow \R$ be an injective function. If $\sigma \neq \sigma'$ are two distinct cells of $X$, then we have $\halo(\sigma;f) \cap \halo(\sigma';f) = \varnothing$.
\end{proposition}
\begin{proof}
	If a cell $\tau$ is contained in the intersection $\halo(\sigma;f) \cap \halo(\sigma';f)$, then by definition both $\sigma$ and $\sigma'$ must lie in $\st^-(\tau;f)$ and attain the minimal value of $f$ over this lower star. This violates the injectivity of $f$.
\end{proof}

The halo of a cell $\sigma$ along an injective $f$ is also well-behaved as a subposet of $\cl(\sigma)$.

\begin{proposition}\label{prop_haloconvex}
Assume $f:X \hookrightarrow \R$ is injective. For each cell $\sigma$ of $X$, the augmented halo $\ahalo(\sigma;f)$ is an up-set of $\cl(\sigma)$.
\end{proposition}
\begin{proof}
We seek to show that given two cells $\beta \geq \alpha$ in $\cl(\sigma)$ with $\alpha$ contained in $\ahalo(\sigma;f)$, we also have $\beta \in \ahalo(\sigma;f)$. If any of the inequalities in $\sigma \geq \beta \geq \alpha$ is an equality, then the result is immediate. Therefore, we may safely assume $\sigma > \beta > \alpha$. We now claim that $f(\beta) > f(\sigma)$; to see why, note by injectivity of $f$ that we either have
\begin{itemize}
	\item $f(\beta) > f(\alpha)$, in which case $f(\beta) > f(\sigma)$ holds because $\sigma \in \st^-(\alpha;f)$; or, 
	\item $f(\beta) < f(\alpha)$, whence $\beta \in \st^-(\alpha;f)$ and so $f(\beta) > f(\sigma)$ by definition of $\argmin$. 
\end{itemize} 
Thus, in both cases we obtain $\sigma \in \st^-(\beta;f)$. Assume, for the sake of contradiction, that $\beta \notin \halo(\sigma;f)$; then there must exist another cell $\sigma' \neq \sigma$ in $\st^-(\beta;f)$ satisfying $f(\sigma') < f(\sigma)$. But any such $\sigma'$ would automatically lie in the lower star of $\alpha$ and contradict the fact that $\sigma$ has the smallest $f$-value amongst lower-star cofaces of $\alpha$. This gives the desired containment $\beta \in \halo(\sigma;f)$ and establishes upward closure of the augmented halo. \end{proof}

Since the removal of any up-set of cells from a regular CW complex leaves behind a subcomplex, the preceding result guarantees that
\begin{align}\label{eq:res}
\res(\sigma;f) := \cl(\sigma)-\ahalo(\sigma;f),
\end{align} is a (possibly empty) subcomplex of $\cl(\sigma)$, and hence, of $X$ whenever $f$ is injective. We will call this subcomplex the {\bf shadow} of $\sigma$ along $f$. Define, for each real number $c$, the sub- and super-levelset of $f$ at $c$ as 
\begin{align*}
	f_{\leq c} &:= \set{\sigma \in X \mid f(\sigma) \leq c}, \\
	f_{\geq c} &:= \set{\sigma \in X \mid f(\sigma) \geq c}.
\end{align*}
Our interest in the halo stems from the fact that it completely determines the differences between closures of certain sublevelsets.

\begin{lemma}\label{lem:halodiff}
	Let $f: X \hookrightarrow \R$ be an injective function, and let $[c,d]\subset\R$ be an interval such that $f_{\geq c} \cap f_{\leq d}$ consists of a single cell $\sigma \in X$. For any such interval, the difference 
	\[
	\Delta := \cl(f_{\leq d}) - \cl(f_{\leq c})
	\] of sublevelsets has the following structure:  
	\begin{enumerate}
		\item if $\sigma \in \cl(f_{\le c})$, then $\Delta = \varnothing$;
		\item if $\sigma \notin \cl(f_{\le c})$, then $\Delta = \ahalo(\sigma;f)$.
	\end{enumerate}
\end{lemma}
\begin{proof}
A cell $\gamma$ lies outside $\cl(f_{\leq c})$ if and only if $f(\tau) > c$ holds for every coface $\tau \geq \gamma$. Conversely, $\gamma$ lies inside $\cl(f_{\leq d})$ if and only if there exists some coface $\tau_* \geq \gamma$ with $f(\tau_*) \leq d$. 
Thus, in order for $\gamma$ to lie in $\Delta$, two properties must hold simultaneously.
First, there must exist some $\tau_* \geq \gamma$ with $f(\tau_*) \in (c,d]$. And second, for every other $\tau \neq \tau_*$ in $\st(\gamma)$ we must have $f(\tau) > c$. The hypothesis requiring $f_{\geq c} \cap f_{\leq d}$ to equal $\set{\sigma}$ affects both of these properties and yields the desired conclusions. Explicitly,
\begin{enumerate} 
	\item The first property simplifies to the requirement that $\sigma = \tau_*$ lies in $\st(\gamma)$ and satisfies $f(\sigma) \neq c$. Now if $\sigma \in \cl(f_{\leq c})$ then all of its faces also lie in $\cl(f_{\leq c})$ by the subcomplex property. But since $\Delta$ can only contain faces of $\sigma$, it must be empty in this case.
	\item The second property strengthens to $f(\tau) > d$ for all $\tau \neq \sigma$ in $\st(\gamma)$. Thus, a cell $\gamma$ lies in $\Delta$ if and only if $\set{\sigma} = \argmin[f/\st(\gamma)]$. There are now two cases to consider. Either $\gamma = \sigma$, or $\gamma < \sigma$ and $f(\gamma) > d \geq f(\sigma)$. In the latter case, note that $\sigma$ lies in $\st^-(\gamma;f)$ and is the argmin of $f$ over this lower star.
\end{enumerate}
To conclude the proof, we note from the strengthened version of the second property above that $\gamma \in \Delta$ holds if and only if $\gamma$ lies in $ \set{\sigma} \sqcup \halo(\sigma;f)$, as desired.
\end{proof}

\section{Stratified Discrete Morse Functions}\label{sec:stratmorse}

Fix a finite regular CW complex $X$ and a function $f:X \to \R$. The {\em upper closure} of each cell $\sigma$ along $f$ is defined as:
\begin{align}\label{eq:upclos}
\cl^+(\sigma;f) := \set{\tau \in X \mid \tau < \sigma \text{ and }f(\tau) \geq f(\sigma)}.
\end{align}
Upper closures are dual to lower stars from \eqref{eq:lowstar} in the sense that $\tau \in \cl^+(\sigma;f)$ holds if and only if $\sigma \in \st^-(\tau;f)$ holds. We call $f$ a {\em discrete Morse function} on $X$ in the sense of Forman  whenever the inequality
\[
\card{\cl^+(\sigma;f)} + \card{\st^-(\sigma;f)} \leq 1
\]
holds for every cell $\sigma \in X$, with $\card{~\cdot~}$ denoting cardinality --- see Definition 2.2 and Lemma 2.5 of \cite{forman1998} for details. Since $X$ is finite, we can always perturb $f$ slightly so that it becomes injective while preserving all relevant cardinalities. Let us now impose a stratification $X_\bullet$ on $X$, as in Section \ref{ssec:strat}. 

\begin{definition}\label{def:sdMf}
A \define{stratified discrete Morse function} on $X$ is any injective map $f:X \hookrightarrow \R$ such that for every cell $\sigma$ lying in a stratum $S \subset X$, we have 
\begin{enumerate}
\item $\card{\cl^+(\sigma;f) \cap S} + \card{\st^-(\sigma;f) \cap S} \leq 1$, and
\item if $\cl^+(\sigma;f)\cap S = \{
\tau\}$, then $\tau\in \halo(\sigma;f)$ and $\cl(\sigma) \scoll \res(\sigma;f)$.
\end{enumerate}
\end{definition}

The first condition of this definition requires the restriction $f|_S:S \to \R$ to be a standard (albeit injective) discrete Morse function on each stratum $S$ in the sense of \cite{forman1998}. The second condition is new to the best of our knowledge. Whenever a cell $\sigma \in S$ admits a unique $\tau$ in $\cl^+(\sigma;f) \cap S$, this condition requires the closure of $\sigma$ to collapse onto the shadow subcomplex from \eqref{eq:res}. In other words, all cells in the augmented halo $\ahalo(\sigma;f)$ must disappear along a sequence of elementary collapses\footnote{Moreover, the first pair in any such sequence is necessarily $(\tau,\sigma)$, because $\sigma$ lies in the open star of every cell in $\cl(\sigma)$ and the only other cell in $S$ that it can possibly be paired with $\sigma$ is $\tau$.}. Another novel feature engendered by our definition is that critical cells come in two flavours. 

\begin{definition}\label{def:crit}
Let $f: X \hookrightarrow \R$ be an injective function satisfying Definition~\ref{def:sdMf}.1. 
A cell $\sigma$ lying in a stratum $S \subset X$ is called
\begin{itemize}
	\item {\bf critical} for $f$ if both ${\cl^+(\sigma;f) \cap S}$ and ${\st^-(\sigma;f) \cap S}$ are empty; and,
	\item {\bf s-critical} for $f$ if it is critical in the sense above and $\st^-(\sigma;f)$ is empty.    
\end{itemize}
(The requirement that $\st^-(\sigma;f) = \varnothing$ is equivalent to $\sigma \notin \halo(\tau;f)$ for any cell $\tau \in X$).
\end{definition}
From \cite[Definition 2.2]{forman1998} it follows that $\sigma$ is critical if and only if it is critical in Forman's sense for the restriction $f|_S:S \to \R$. If the stratification is trivial, then the only stratum in sight is $S = X$; in this case, s-criticality coincides with criticality because the first requirement of the above definition already forces $\st^-(\sigma;f)$ to be empty.

\section{Local Morse Data} \label{sec:morselems}

We now turn towards the task of establishing the fundamental results pertaining to sublevelsets of a stratified discrete Morse function $f:X \hookrightarrow \R$. The first of these pertains to the case where an interval $[c,d] \subset \R$ contains no $s$-critical values. The goal, as mentioned in the preamble to Lemma \ref{lem:halodiff}, is to relate  $\cl(f_{\leq c})$ and $\cl(f_{\leq d})$.

\begin{theorem}\label{thm:a}
Let $f: X \hookrightarrow \R$ be a stratified discrete Morse function and consider an interval $[c,d] \subset \R$ for which $f_{\geq c} \cap f_{\leq d}$ consists of a single cell $\sigma$, which 
is not s-critical.
Then there exists a filtered collapse $\cl(f_{\le d}) \scoll~ \cl(f_{\le c})$.
\end{theorem}
\begin{proof}
If $f(\sigma) = c$ then by the uniqueness of $\sigma$ we get $f_{\leq c} = f_{\leq d}$ and there is nothing to prove; thus we assume $f(\sigma) \in (c,d]$. The argument now decomposes into two cases: either the lower star $\st^-(\sigma;f)$ is empty, or it is not. Of these, the nonempty case is far simpler.

{\bf Case 1:} If $\st^-(\sigma;f)$ contains some cell $\tau$, then by \eqref{eq:lowstar} we have both $\sigma < \tau$ and $f(\sigma) \geq f(\tau)$. By the injectivity of $f$ and the uniqueness of $\sigma$ in $f_{\geq c} \cap f_{\leq d}$, we get $f(\tau) < c$, which forces $\sigma \in \cl(f_{\leq c})$. Lemma \ref{lem:halodiff}(1) now gives $\cl(f_{\le d})=\cl(f_{\le c})$, so the trivial collapse suffices. 

{\bf Case 2:} If $\st^-(\sigma;f)$ is empty, then since $f(\sigma) > c$ we have $\sigma \notin \cl(f_{\leq c})$. Therefore Lemma \ref{lem:halodiff}(2) yields
\[
\cl(f_{\leq d}) - \cl(f_{\leq c}) = \ahalo(\sigma;f).
\]
Since $\sigma$ is not s-critical and $\st^-(\sigma;f)$ is empty, Definition \ref{def:sdMf}(2) forces $\cl(\sigma) \scoll \res(\sigma;f)$. We now claim that the same sequence of elementary collapses which reduces $\cl(\sigma)$ to $\res(\sigma;f)$ also reduces $\cl(f_{\leq d})$ to $\cl(f_{\leq c})$. First, note from \eqref{eq:res} that
\[
\cl(\sigma) - \res(\sigma;f) =  \ahalo(\sigma;f),
\] so in both scenarios one must excise only those cells which lie in the augmented halo of $\sigma$. It remains to confirm that the free-face pairs removed from $\cl(\sigma)$ are also free-face pairs in $\cl(f_{\leq d})$. This amounts to checking open stars, as described in Section \ref{ssec:collapse}. Explicitly, for each cell $\gamma \in \ahalo(\sigma;f)$, we must establish
\[
\big[\st(\gamma) \cap \cl(\sigma)\big] \supset \big[\st(\gamma) \cap \cl(f_{\leq d})\big],
\] 
since the reverse containment is guaranteed by the fact that $f(\sigma) \leq d$. We will proceed by contradiction, assuming the existence of some $\rho > \gamma$ lying in $\cl(f_{\leq d}) - \cl(\sigma)$. By definition of $\cl(f_{\leq d})$ there exists a coface $\rho' \geq \rho$ with $f(\rho') \leq d$. Now $\rho \notin \cl(\sigma)$ implies $\rho' \neq \sigma$; and since $\sigma$ is the only cell valued in $[c,d]$, we get $f(\rho') < c \leq f(\sigma)$. This violates the argmin property required by $\gamma \in \ahalo(\sigma;f)$ and hence concludes the proof.
\end{proof}

One strange consequence of the above result is that critical cells of $f$ can only affect the filtered topology of sublevelsets if they also happen to be s-critical. We now quantify the change in topology from $\cl(f_{\leq c})$ to $\cl(f_{\leq d})$ whenever the interior $(c,d)$ contains a single s-critical value. This is the purview of the second fundamental theorem of Morse theory. Since our $f$ is injective, we may safely restrict attention to the case where $f_{\geq c} \cap f_{\leq d}$ is a singleton.  

\begin{theorem}\label{thm:b}
Let $f: K \hookrightarrow \R$ be a stratified discrete Morse function and $\sigma \in X$ an s-critical cell with $c := f(\sigma)$. For all sufficiently small $\eps > 0$, the subcomplex $\cl(f_{\leq c+\eps})$ is the union of $\cl(\sigma)$ with $\cl(f_{\leq c-\eps})$ along the intersection $\res(\sigma;f)$.
\end{theorem}
\begin{proof}
By finiteness of $X$ and injectivity of $f$, we are free to choose $\eps > 0$ such that $\sigma$ is the only cell of $X$ lying in $f_{\geq c-\eps} \cap f_{\leq c+\eps}$. And since $\sigma$ is s-critical, we know from Definition \ref{def:crit} that $\st^-(\sigma;f)$ is empty. Consequently, $\sigma$ is not in $\cl(f_{\le c-\eps})$, so Lemma \ref{lem:halodiff}(2) implies that $\cl(f_{\leq c+\eps})$ equals the union $\cl(f_{\leq c-\eps}) \cup \ahalo(\sigma;f)$. Since $\ahalo(\sigma;f) \subset \cl(\sigma)$ holds by Definition \ref{def:halo}, we get the containment 
\[\cl(f_{\leq c+\eps}) \subset \big[\cl(f_{\leq c-\eps}) \cup \cl(\sigma)\big].
\] Conversely, $\cl(\sigma) \subset \cl(f_{\leq c+\eps})$ since $f(\sigma) = c \leq c+\eps$, so the reverse containment also holds. This proves that $\cl(f_{\leq c+\eps})$ equals the union displayed above. It remains to establish the pushout property by showing that
\[
\res(\sigma;f) = \cl(f_{\leq c-\eps}) \cap \cl(\sigma).
\] Below we will show that the left side is contained in the right side and vice-versa. 

$\star~{\text{\bf left} \subset \text{\bf right}}$ : By \eqref{eq:res} and Definition \ref{def:halo}, if a cell $\gamma$ lies in $\res(\sigma;f)$ then $\gamma < \sigma$ and 
\[
\sigma \neq \argmin[f/\st^-(\gamma;f)].
\] 
We claim that this inequality forces $\gamma \in \cl(f_{\leq c- \eps})$. 
If $\st^-(\gamma;f)$ is empty, then $f(\gamma) < f(\sigma)$ and hence $f(\gamma) \leq c-\eps$ and we get the desired $\gamma \in f_{\leq c-\eps} \subset \cl(f_{\leq c-\eps})$. Alternately, if some $\rho \in \st^-(\gamma;f)$ satisfies $f(\rho) < f(\sigma)$, then $f(\rho) \leq c-\eps$ and hence $\gamma \in \cl(f_{\leq c-\eps})$.  

$\star~{\text{\bf right} \subset \text{\bf left}}$: We show that any cell $\gamma$ in $\cl(\sigma) \cap \cl(f_{\leq c-\eps})$ must lie outside $\ahalo(\sigma;f)$. If this $\gamma$ lies in $f_{\leq c-\eps}$ then we have $\sigma \notin \st^-(\gamma;f)$ since $\sigma$ has a higher $f$-value. Thus, $\sigma$ can not be the argmin of $f$ over the lower star of $\gamma$, which means that $\gamma$ lies outside the augmented halo of $\sigma$. On the other hand, if $\gamma$ lies in $\cl(f_{\leq c-\eps})$ because some cell $\rho > \gamma$ satisfies $f(\rho) \leq c-\eps$, then the presence of this $\rho$ once again rules out the possibility that $\sigma = \argmin[f/\st^-(\sigma;f)]$ because we have 
\[
f(\rho) \leq c-\eps < c = f(\sigma).
\] Consequently, we see that $\gamma$ must lie in $\cl(\sigma) - \ahalo(\sigma;f)$ and hence in $\res(\sigma;f)$ by \eqref{eq:res}.
\end{proof}

\section{Local Morse Data via the Upper Envelope}\label{sec:envelope}

The {\bf barycentric subdivision} of a finite regular CW complex $X$ is the simplicial complex $\sd(X)$ whose $k$-simplices, for all $k \geq 0$, are all strictly ascending chains of the form 
\[
\xi := [\sigma_0 < \sigma_1 < \cdots < \sigma_k]
\]
in the face poset $(X,\leq)$. The faces $\eta < \xi$ are obtained by deleting one or more cells from the underlying chain $\sigma_0 < \cdots < \sigma_k$. There is a natural order-preserving {\em last cell} map $\omega:\sd(X) \twoheadrightarrow X$ that sends each such $\xi$ to its terminal cell $\sigma_k$. It is relatively straightforward to check that $\sd(X)$ canonically inherits a stratification from $X$ via $\omega$ --- namely, for each stratum $S \subset X$ there exists a unique stratum $S' \subset \sd(X)$ such that the simplex $\xi$ lies in $S'$ if and only if its last cell $\omega(\xi)$ lies in $S$. It is readily seen that $S'$ consists exclusively of barycentrically subdivided $S$-cells, so we will denote it by $\sd(S)$ rather than $S'$. We assume throughout this section that $X$ -- and hence, $\sd(X)$ -- is stratified.
\begin{definition}\label{def:env}
The {\bf upper envelope} of an injective function $f:X \hookrightarrow \R$ is the function $\U{f}:\sd(X) \to \R$ given by
\[
\U{f}\left([\sigma_0<\cdots<\sigma_k]\right) := \max\set{f(\sigma_i)\mid 0 \leq i \leq k}.
\]
(One may safely define upper envelopes for arbitrary $X \to \R$, but the only case of interest to us here occurs when $f$ is injective).
\end{definition}
 The reader is warned that $\U{f}$ is almost never injective even though the original $f$ is injective. Consequently, upper envelopes of stratified discrete Morse functions on $X$ rarely yield stratified discrete Morse functions on $\sd(X)$. The advantage of passing to the upper envelope is that $\U{f}$ is order-preserving --- given $\xi \leq \eta$ in $\sd(X)$, we always have $\U{f}(\xi) \leq \U{f}(\eta)$ in $\R$. This monotonicity fuels the following elementary result. 

\begin{proposition}\label{prop:sdsubcomp}
Given any function $f:X \to \R$ and value $c \in \R$, the sublevelset
\[
\U{f}_{\leq c} := \set{\xi \in \sd(X) \mid \U{f}(\xi) \leq c}
\]
is a subcomplex of $\sd(X)$.
\end{proposition}
\begin{proof}
Assume a simplex $\xi \in \sd(X)$ lies in $\U{f}_{\leq c}$ and consider a face $\eta < \xi$. We see from Definition \ref{def:env} that $\U{f}(\eta)$ is the maximum of $f$ over a smaller set, and therefore we obtain $\U{f}(\eta) \leq \U{f}(\xi) \leq c$. Therefore, we have $\eta \in \U{f}_{\leq c}$ as desired. 
\end{proof}

We now seek to explicitly describe the local Morse data of a stratified Morse function $f:X \hookrightarrow \R$ at a cell $\sigma \in X$ in terms of the sublevesets of $\U{f}$ in $\sd(X)$ rather than (closures of) sublevelsets of $f$ in $X$. For this purpose, it is convenient to better understand the neighbourhood of the corresponding vertex $[\sigma] \in \sd(X)$. The open star of $[\sigma]$ is evidently the collection of all barycentric simplices whose underlying chains contain $\sigma$ in some position:
\[
\st[\sigma] = \set{[\sigma_0 < \cdots < \sigma_k] \in \sd(X) \mid \sigma_i = \sigma \text{ for some }0 \leq i \leq k}.
\]
It follows that the {\bf link} of $[\sigma]$, given by
\[
\lk[\sigma] := \cl(\st[\sigma]) - \st[\sigma],
\] is the subcomplex of $\sd(X)$ comprising all $\xi$ whose underlying chain $\sigma_0 < \cdots < \sigma_k$ can be augmented by $\sigma$. 

\begin{remark} \label{rem:linktype} Specifically, given any $\xi = [\sigma_0 < \cdots < \sigma_k]$ in $\lk[\sigma]$, we must have either
\begin{enumerate}
    \item $\sigma < \sigma_0$, so $\sigma$ fits at the beginning of the chain; or, 
    \item $\sigma_k < \sigma$, so $\sigma$ fits at the end of the chain; or,
    \item $\sigma_i < \sigma < \sigma_{i+1}$ for some $i$ in $\set{0,1,\ldots,k-1}$.
\end{enumerate}
In every case these inequalities are strict, so $\sigma \neq \sigma_i$ for any $i$ and therefore $\xi \notin \st[\sigma]$. 
\end{remark} 
To confirm that $\lk[\sigma]$ is a subcomplex of $\sd(X)$ as asserted above, one can either use the fact that $\st[\sigma]$ is an up-set in its closure, or simply note that the removal of any $\sigma_i$ from $\sigma_0 < \cdots < \sigma_k$ preserves $\sigma$-augmentability.  
Set $c := f(\sigma)$, and define the {\bf lower link} of $[\sigma]$ along $\U f$ as 
\begin{align} \label{eq:lowlink}
\lk^-_{{f}}[\sigma] := \lk[\sigma] \cap \U{f}_{< c},
\end{align}
where $\U{f}_{< c}$ is defined analogously to $\U{f}_{\leq c}$. Since $X$ is finite, there exists $\eps>0 $ such that $\U f_{<c}$ coincides with $\U f_{\le c-\eps}$, and hence forms a subcomplex of $\sd(X)$ by Proposition \ref{prop:sdsubcomp}. It follows that the lower link is an intersection of subcomplexes, and hence a subcomplex of $\sd(X)$. Lower links play a crucial role in the study of sublevelsets of injective real-valued maps defined on the vertices of a simplicial complex. 
By far most prominent instance of this phenomenon is found in the piecewise linear Morse theory of Bestvina and Brady \cite{bb}. 

\begin{definition}\label{def:join}
Let us denote by $K_i$ the set of $i$-dimensional simplices in a simplicial complex $K$. Given another simplicial complex $L$, we recall that the {\bf join} $K \star L$ is the new complex whose $n$-simplices are given by
\[
(K \star L)_n = K_n \sqcup L_n \sqcup \coprod_{i+j=n-1} K_i \times L_{j}.
\] 
The {\bf cone} over $L$ is the special case of this construction where $K$ consists of a single vertex.
\end{definition}

The cone $[\sigma] \star \lk_f^-[\sigma]$ is evidently a subcomplex of $\sd(X)$; combined with the lower link itself, it completely characterises the inclusion of sublevelsets $\U{f}_{\leq c-\eps} \hookrightarrow \U{f}_{\leq c+\eps}$ when $f$ is injective on $X$.

\begin{proposition}\label{prop:sdpushout}
Let ${\U{f}: \sd(X)\to \R}$ be the upper envelope of an injective map $f:X \hookrightarrow \R$. For every cell $\sigma \in X$ with $c:= f(\sigma)$, there exists sufficiently small $\eps > 0$ such that  $\U{f}_{\leq c+\eps}$ is the union of $\U{f}_{\leq c-\eps}$ with $[\sigma] \star \lk^-_f[\sigma]$ along their intersection $\lk^-_f[\sigma]$.
\end{proposition}
\begin{proof}
Since $f$ is injective and $X$ is finite, we may safely choose an $\eps > 0$ such that $\sigma$ is the unique cell of $X$ lying in the intersection $f_{\geq c-\eps} \cap f_{\leq c+\eps}$. By Definition \ref{def:env}, the upper envelope $\U{f}$ attains the same values as $f$ in $\R$. Therefore, any cell 
\[
\eta = [\sigma_0 < \cdots < \sigma_k] ~ \text{ lying in } \left(\U{f}_{\leq c+\eps} - \U{f}_{\leq c-\eps} \right)
\] satisfies $\U{f}(\eta) = c$. It now follows from Definition \ref{def:env} and the injectivity of $f$ that  $\sigma = \sigma_i$ for some $0 \leq i \leq k$ and $f(\sigma_j) < c$ for all the other $j \neq i$. Therefore, using the uniqueness of $\sigma$ guaranteed by our choice of $\eps$, we know that the face $\xi < \eta$ obtained by removing $\sigma_i$ from the underlying chain satisfies $f(\xi) \leq c-\eps$, and hence $\xi \in \lk^-_f[\sigma]$. Thus, our $\eta$ has the form $\xi \cup [\sigma]$ for some $\xi \in \lk^-_f[\sigma]$, as desired.
\end{proof}

The preceding result may disappoint the reader who seeks a more faithful discretisation of $f$'s local Morse data   at $[\sigma]$ \`{a} la Goresky and MacPherson \cite{GoreMPer1988}. The central discrepancy here is that the lower link $L$ of $[\sigma]$ along $f$ does not admit an obvious factorisation of the form $TL \times NL$, where $TL$ is a {\em tangential} component which lies entirely within the unique stratum $\sd(S) \subset X$ containing $[\sigma]$, while $NL$ is the {\em normal} component that resides in strictly higher strata $\sd(T) > \sd(S)$. There are two immediate obstacles --- first, simplicial complexes are not closed under products, so it is unreasonable to expect the desired $TL$ and $NL$ to form subcomplexes of $\sd(X)$. Secondly, the simplicial setting frustrates all na\"ive attempts to take arbitrarily small neighbourhoods around $[\sigma]$ --- the smallest available option is $\cl(\st[\sigma])$, which will necessarily intersect  {\em lower} strata $\sd(T) < \sd(S)$ whenever there exists a face $\tau < \sigma$ which lies in $T$.

Nevertheless, we persist. Consider a (not necessarily critical) cell $\sigma$ of $X$ and fix 
an injective function $f:X \hookrightarrow \R$. For brevity, we write $L$ for the lower link $\lk^-_f[\sigma]$ from \eqref{eq:lowlink}. The first step is to isolate two relevant subcomplexes of $L$, which correspond to the first two cases of Remark \ref{rem:linktype}.
\begin{definition}\label{def:HV}
A simplex $\xi = [\sigma_0 < \cdots < \sigma_k]$ in $L$ is said to lie in   
\begin{enumerate}
\item the {\bf horizontal part} $H$ of $L$ whenever $\sigma_k < \sigma$, and in
\item the {\bf vertical part} $V$ of $L$ whenever $\sigma < \sigma_0$.
\end{enumerate}
(Both defining properties are preserved when passing to faces of $\xi$, whence $H$ and $V$ form subcomplexes of $L$). 
\end{definition}

We catalogue the following elementary and helpful properties of $H$ and $V$. 
\begin{proposition}
    Let $H$ and $V$ be the horizontal and vertical parts of the lower link $L$ of $\sigma$ along $f$. Then,
    \begin{enumerate}
    \item the intersection $H \cap V$ is empty;
    \item the join $H \star V$ is $L$.
    \end{enumerate}
\end{proposition}
\begin{proof} For (1), note that if a simplex $\xi = [\sigma_0 < \cdots < \sigma_k]$ lies in $H \cap V$ then we arrive at the absurdity $\sigma < \sigma_0 < \sigma_k < \sigma$, where the first inequality comes from the definition of $V$ and the last inequality comes from the definition of $H$. Moving to (2), assume that $\xi$ lies in the difference $L - (H \cup V)$. By Remark \ref{rem:linktype} there is some $i$ in $\set{0,1,\ldots,k-1}$ with $\sigma_i < \sigma < \sigma_{i+1}$. But now, the front face $[\sigma_0 < \cdots < \sigma_i]$ of $\xi$ lies in $H$ while the back face $(\sigma_{i+1}< \cdots < \sigma_k)$ lies in $V$, so $\xi$ is indeed the disjoint union of simplex of $H$ and a simplex of $V$ as required by Definition \ref{def:join}.
\end{proof}

Combining the preceding result with Proposition \ref{prop:sdpushout}, we see that the local Morse data for $\U{f}$ at $\sigma$ is $\big(H \star V\big) \hookrightarrow \big([\sigma] \star (H \star V)\big)$. Since the join operation is manifestly associative on simplicial complexes, we may safely rewrite this as
\begin{align}\label{eq:HVmorsdata}
\big(H \star V\big) \hookrightarrow \big(([\sigma] \star H) \star V\big).
\end{align}
Our next task is to examine how $H$ and $V$ interact with the stratification of $\sd(X)$ inherited from $X$. From here onwards, we will assume that
$f: X\hookrightarrow\R$ satisfies Definition~\ref{def:sdMf} and that
$\sigma$ is a critical cell of $f$, but do not insist on s-criticality.
\begin{proposition} \label{prop:HVstrata}
Let $H$ and $V$ be the horizontal and vertical part of the link $L$ of $[\sigma] \in \sd(X)$. If $S \subset X$ is the stratum containing $\sigma$, then
    \begin{enumerate}
     \item $H$ lies in the union $\bigcup_{T \leq S} \sd(T)$ of (subdivided) strata lying at or below the level of $S$ in the frontier partial order, while
     \item $V$ lies in the union $\bigcup_{T > S} \sd(T)$ of (subdivided) strata lying strictly above $S$ with respect to the frontier partial order.
    \end{enumerate}
    (In fact the first assertion holds for all $\sigma$, whereas the second assertion requires $\sigma$ to be critical). 
\end{proposition}
\begin{proof}
To establish (1), we note from Definition \ref{def:HV} that a simplex $\xi = [\sigma_0 < \cdots < \sigma_k]$ of $L$ lies in $H$ if and only if the last cell $\omega(\xi) = \sigma_k$ is a strict face of $\sigma$ in $X$.
Thus if $T \subset X$ is the stratum containing $\sigma_k$, we get a nonempty $\cl(S) \cap T$ and hence $T \leq S$. Passing to the subdivision, $\xi$ must lie in $\sd(T)$ as desired. For (2), given $\xi \in V$, we know that $\sigma$ is a strict face of the first cell $\sigma_0$ by definition of $V$, and that $f(\sigma_0) < f(\sigma)$ holds because $\xi$ lies in the lower link $L$ of $\sigma$. Therefore, we have $\sigma_0 \in \st^-(\sigma;f)$, and the fact that $\sigma$ is critical now forces $\sigma_0 \notin S$ (as per Definition \ref{def:crit}). Thus, $\sigma_0$ lies in some $T_0 > S$ along the frontier order. Proceeding similarly from left to right along the chain $\sigma_0 < \cdots < \sigma_k = \omega(\xi)$, we see that the stratum $T_i$ containing $\sigma_i$ must satisfy $T_i \geq T_{i-1}$ for all $i > 0$. In particular, it follows that the stratum $T_k$ whose subdivision contains $\xi$ satisfies $T_k \geq T_0 > S$, as desired.
\end{proof}

Finally, we combine \eqref{eq:HVmorsdata} with Propositions \ref{prop:sdpushout}  and \ref{prop:HVstrata} to arrive at the main result of this Section.

\begin{theorem}\label{thm:tndecomp}
Let $\sigma$ be a critical cell for an injective function $f:X \hookrightarrow \R$ which satisfies the first requirement of Definition \ref{def:sdMf}. Set $c:=f(\sigma)$ and choose any $\eps>0$ so that 
$\sigma$ is the unique cell of $X$ with $f(\sigma)\in[c-\eps,c+\eps]$. Let $H$ and $V$ be the horizontal and vertical parts of the lower link of $[\sigma]$ along $\U{f}$. Then,
\begin{enumerate}
    \item $\U{f}_{\leq c+\eps}$ is the union of $\U{f}_{\leq c-\eps}$ with $([\sigma] \star H) \star V$ along their intersection $H \star V$; and,
    \item if $S \subset X$ is the stratum containing $\sigma$, then $[\sigma] \star H$ lies in the union of strata $\leq \sd(S)$ while $V$ lies in the union of strata $> \sd(S)$.
\end{enumerate}
\end{theorem}
\begin{proof}
The only assertion which has not been explicitly addressed in previous results is the claim that $[\sigma] \star H$ lies in the union $U$ of strata $\leq \sd(S)$. We already know that every simplex $\xi \in H$ lies in $U$ from Proposition \ref{prop:HVstrata}, so it suffices to check that $[\sigma] \cup \xi$ is also in $U$. By Definition \ref{def:HV}, we see that $\sigma$ must be the last cell in the underlying chain of $[\sigma] \cup \xi$, whence $[\sigma] \cup \xi$ always lies in $S$ and hence in $U$, as desired. 
\end{proof}

\bibliographystyle{abbrv}
\bibliography{bibliography}

\begin{thebibliography}{10}

\bibitem{adiprasito}
K.~Adiprasito and B.~Benedetti.
\newblock Collapsibility of {CAT}(0) spaces.
\newblock {\em Geometriae Dedicata}, 206:181--199, 2020.

\bibitem{atiyahbott}
M.~F. Atiyah and R.~Bott.
\newblock The {Y}ang--{M}ills equations over {R}iemann surfaces.
\newblock {\em Philosophical Transactions of the Royal Society of London. Series A}, 308(1505):523--615, 1983.

\bibitem{babson}
E.~Babson and P.~Hersh.
\newblock Discrete {M}orse functions from lexicographic orders.
\newblock {\em Transactions of the American Mathematical Society}, 357(2):509--534, 2005.

\bibitem{bb}
M.~Bestvina and N.~Brady.
\newblock Morse theory and finiteness properties of groups.
\newblock {\em Inventiones Mathematicae}, 129:445--470, 1997.

\bibitem{chari}
M.~K. Chari.
\newblock On discrete {M}orse functions and combinatorial decompositions.
\newblock {\em Discrete Mathematics}, 217(1):101--113, 2000.

\bibitem{cohen}
M.~M. Cohen.
\newblock {\em A course in simple homotopy theory}.
\newblock Springer-Verlag, 1973.

\bibitem{conley}
C.~C. Conley.
\newblock {\em Isolated Invariant Sets and the Morse Index}, volume~38 of {\em CBMS Regional Conference Series in Mathematics}.
\newblock American Mathematical Society, 1978.

\bibitem{cornea}
O.~Cornea, K.~A. {De Rezende}, and M.~R. {Da Silveira}.
\newblock Spectral sequences in {C}onley’s theory.
\newblock {\em Ergodic Theory and Dynamical Systems}, 30(4):1009–1054, 2010.

\bibitem{dmsheaf}
J.~Curry, R.~Ghrist, and V.~Nanda.
\newblock Discrete {M}orse theory for computing cellular sheaf cohomology.
\newblock {\em Foundations of Computational Mathematics}, 16(4):875--897, 2016.

\bibitem{braid}
D.~Farley and L.~Sabalka.
\newblock Discrete {M}orse theory and graph braid groups.
\newblock {\em Algebraic \& Geometric Topology}, 5:1075--1109, 2005.

\bibitem{floer}
A.~Floer.
\newblock Morse theory for {L}agrangian intersections.
\newblock {\em Journal of Differential Geometry}, 28(3):513--547, 1988.

\bibitem{forman1998}
R.~Forman.
\newblock Morse theory for cell complexes.
\newblock {\em Advances in mathematics}, 134(1):90--145, 1998.

\bibitem{Freij2009}
R.~Freij.
\newblock Equivariant discrete {M}orse theory.
\newblock {\em Discrete Mathematics}, 309(12):3821--3829, 2009.

\bibitem{GoreMPer1988}
M.~Goresky and R.~MacPherson.
\newblock {\em Stratified Morse theory}.
\newblock Springer, 1988.

\bibitem{welker}
M.~J\"{o}llenbeck and V.~Welker.
\newblock {\em Minimal Resolutions via Algebraic Discrete {M}orse Theory}, volume 197 of {\em Memoirs of the American Mathematical Society}.
\newblock AMS, 2009.

\bibitem{knudwang2022}
K.~Knudson and B.~Wang.
\newblock Discrete stratified {M}orse theory: Algorithms and a user’s guide.
\newblock {\em Discrete \& Computational Geometry}, 67(4):1023--1052, 2022.

\bibitem{kukiela}
M.~Kukiela.
\newblock The main theorem of discrete {M}orse theory for {M}orse matchings with finitely many rays.
\newblock {\em Topology and its Applications}, 160(9):1074--1082.

\bibitem{LipKubMroWan201}
M.~Lipiński, J.~Kubica, M.~Mrozek, and T.~Wanner.
\newblock Conley-{M}orse-{F}orman theory for generalized combinatorial multivector fields on finite topological spaces.
\newblock {\em Journal of Applied and Computational Topology}, 7:139--184, 2023.

\bibitem{mather}
J.~N. Mather.
\newblock Notes on topological stability.
\newblock {\em Bulletin of the American Mathematical Society}, 49(4):475--506, 2012.

\bibitem{dmpers}
K.~Mischaikow and V.~Nanda.
\newblock Morse theory for filtrations and efficient computation of persistent homology.
\newblock {\em Discrete \& Computational Geometry}, 50(2):330--353, 2013.

\bibitem{morse}
M.~Morse.
\newblock The calculus of variations in the large.
\newblock {\em Annals of Mathematics}, 26(3):213--234, 1925.

\bibitem{mro2017}
M.~Mrozek.
\newblock Conley–{M}orse–{F}orman theory for combinatorial multivector fields on lefschetz complexes.
\newblock {\em Foundations of Computational Mathematics}, 17:1585--1633, 2017.

\bibitem{MrozekWanner2015Dynamics}
M.~Mrozek and T.~Wanner.
\newblock Dynamics of combinatorial multivector fields.
\newblock In {\em Connection Matrices in Combinatorial Topological Dynamics}, SpringerBriefs in Mathematics, pages 105--118. Springer, Cham, 2015.

\bibitem{flow2}
V.~Nanda.
\newblock Discrete {M}orse theory and localization.
\newblock {\em Journal of Pure and Applied Algebra}, 223(2):459--488, 2019.

\bibitem{locstrat}
V.~Nanda.
\newblock Local cohomology and stratification.
\newblock {\em Foundations of Computational Mathematics}, 20(2):195--222, 2020.

\bibitem{flow1}
V.~Nanda, D.~Tamaki, and K.~Tanaka.
\newblock Discrete {M}orse theory and classifying spaces.
\newblock {\em Advances in Mathematics}, 340:723--790, 2018.

\bibitem{salamon}
D.~Salamon.
\newblock Connected simple systems and the {C}onley index of isolated invariant sets.
\newblock {\em Bulletin of the American Mathematical Society}, 291:1--41, 1985.

\bibitem{salvetti}
M.~Salvetti.
\newblock Topology of the complement of real hyperplanes in ${\mathbb c}^n$.
\newblock {\em Inventiones Mathematicae}, 88(3):603--618, 1987.

\bibitem{weinberger}
S.~Weinberger.
\newblock {\em The Topological Characterization of Stratified Spaces}.
\newblock The University of Chicago Press, 1994.

\bibitem{whitney}
H.~Whitney.
\newblock Tangents to an analytic variety.
\newblock {\em Annals of Mathematics}, 81(3):496--549, 1965.

\bibitem{witten}
E.~Witten.
\newblock Supersymmetry and morse theory.
\newblock {\em Journal of Differential Geometry}, 17(4):661--692, 1982.

\bibitem{dmcog}
N.~Yerolemou and V.~Nanda.
\newblock Morse theory for complexes of groups.
\newblock {\em Journal of Pure and Applied Algebra}, 228(6):107606, 2024.

\end{thebibliography}

\appendix

\section{Multivector Fields and the Conley Index}\label{sec:conley}

Let $X$ be a finite regular CW complex. A \define{multivector} of $X$ is any nonempty, convex and connected subset $M$ of the poset $(X,\leq)$. A partition $X = \coprod_i M_i$ into multivectors is a \define{multivector field}; these have found substantial use in the combinatorial study of dynamical systems \cite{mro2017, MrozekWanner2015Dynamics, LipKubMroWan201}. 
Here are four examples of multivector fields, in approximately ascending order of generality:
\begin{enumerate}
\item The most trivial case occurs whenever each multivector is a single cell of $X$, so the multivector field is simply the given partition of $X$ into its constituent cells.
\item Given a discrete Morse function $f:X \to \R$ in the sense of Forman, each multivector is either a single $f$-critical cell or a pair $\set{\sigma,\tau}$ of cells where $\sigma < \tau$, $\dim \tau=\dim \sigma+1$ and $f(\sigma) \geq f(\tau)$. The resulting multivector field, called an {\em acyclic partial matching}, is fundamental in discrete Morse theory \cite{chari}.
\item In \cite[Sec 3.2]{Freij2009}, Freij considers a partition of $X$ into {\em intervals}; each interval $I_\sigma^\tau$ consists of all cells $\gamma \in X$ of the form $\sigma \leq \gamma \leq \tau$ for some fixed $\sigma \leq \tau$ in $X$, not necessarily of codimension $\leq 1$. Such interval partitions are also multivector fields.
\item Let $X_\bullet$ be a stratification of $X$ as described in Section \ref{ssec:strat} above. The collection of strata constitutes a multivector field on $X$ --- connectedness holds by definition and convexity follows from Proposition \ref{prop:convexstrata}. 
\end{enumerate}
Conversely, not every multivector field produces a stratification in our sense because the frontier axiom may not hold. 

The {\em exit set} of a multivector $M \subset X$ is  $\ex(M) := \cl(M) - M$; since $M$ is convex, it forms an up-set in its closure, whence $\ex(M)$ is always a (possibly empty) subcomplex of $\cl(M)$. The following definition forms the starting point in a discrete version of Conley's celebrated topological study \cite{conley} of isolated invariant sets --- see for instance \cite[Eq (7.4)]{MrozekWanner2015Dynamics}.

\begin{definition}\label{def:conley}
The {\bf homological Conley index} of a multivector $M \subset X$ is defined as the relative homology group
\[
\Con_k(M) := \HG_k(\cl(M),\ex(M);\Z)
\]
for each dimension $k \geq 0$. We say that $M$ is {\bf critical} whenever $\Con_k(M) \neq 0$ for some $k$.
\end{definition}
\noindent It is readily checked (by writing out the relevant chain complexes) that when $M = \set{\sigma}$ is a single cell, then $\Con_k(M)$ is nontrivial iff $k = \dim \sigma$. On the other hand, if $M = \set{\sigma,\tau}$ contains two cells with $\sigma < \tau$ and $\dim \tau=\dim \sigma+1$, then $\Con_k(M)$ is trivial for all $k \geq 0$. 

Consider the binary relation $\blacktriangleright$ on a multivector field $\cM = \set{M_i \subset X}$ given by 
\[
M_i \blacktriangleright M_j \text{ if there exist } \sigma_i \in M_i \text{ and } \sigma_j \in M_j \text{ with } \sigma_i \geq \sigma_j,
\]
(This is clearly reflexive.) We say that $\cM$ is {\bf acyclic} whenever the transitive closure $\rhd$ of $\blacktriangleright$ forms a partial order on the set of multivectors. Acyclic multivector fields therefore generalise the acyclic partial matchings induced by discrete Morse functions \cite{chari}. In fact, it is easily checked that acyclic multivector fields also subsume the stratifications from Section \ref{ssec:strat}.

\begin{proposition}\label{prop:stratamf}
The strata in any stratification of $X$ form an acyclic multivector field, and $\rhd$ coincides with the frontier partial order on strata.
\end{proposition}
\begin{proof}
By definition, $S \blacktriangleright T$ holds for a pair of strata $S \neq T$ if and only if there exist cells $\sigma \in S$ and $\tau \in T$ with $\sigma \geq \tau$, which in turn occurs if and only if $\cl(S) \cap T$ is nonempty. The desired conclusion now follows from Proposition \ref{prop:frontier}.
\end{proof}

The reader is warned that the preceding result does not imply that the frontier axiom holds amongst the multivectors of an acyclic multivector field; we only recover the partial order implied by the frontier axiom via Proposition \ref{prop:frontier}. On the other hand, acyclicity does allow us to recover the global homology of $X$ from the Conley indices of its multivectors. Explicitly, since $\rhd$ is a partial order, we may enumerate the multivectors $\set{M_1,\ldots,M_n}$ such that $M_i \rhd M_j$ forces $i \geq j$. Therefore, 
\[
F_q := \bigcup_{i \leq q} M_i
\]
constitutes a filtration of $X$ by subcomplexes. We may now examine the associated {\bf Conley-Morse spectral sequence} of this filtration \cite{salamon, cornea}, which is populated by relative homology groups
\[
E^1_{p,q} := \HG_{p+q}(F_q,F_{q-1};\Z),
\]
and converges to $E^\infty_{p,q} \simeq \HG_{p+q}(X)$. Crucially, it follows from excision that $E^1_{p,q}$ is isomorphic to the Conley index $\Con_{p+q}(M_q)$. Proposition \ref{prop:stratamf} guarantees a similar mechanism for assembling the global homology of $X$ from the Conley indices of its strata.

\end{document}